\documentclass[12pt]{amsart}
\usepackage{amsmath,amsthm,latexsym,amscd,amsbsy,amssymb}
 \relax

\chardef\bslash=`\\ 

\makeatletter
\def\verbatim{\interlinepenalty\@M \@verbatim
  \leftskip\@totalleftmargin\advance\leftskip2pc
  \frenchspacing\@vobeyspaces \@xverbatim}
\makeatother
\makeatletter
  \def\dgt@k{\dg@DX=-3 \dg@DY=2 \dg@SIZE=3}
\makeatother

\makeatletter
  \def\dgt@kk{\dg@DX=3 \dg@DY=-1 \dg@SIZE=3}%
\makeatother

\theoremstyle{plain}
\newtheorem{thm}{Theorem}[section]

\newtheorem{lem}[thm]{Lemma}
\newtheorem{pro}[thm]{Proposition}

\theoremstyle{definition}
\newtheorem{rem}[thm]{Remark}
\newtheorem{defin}[thm]{Definition}

\newtheorem*{2loc}{Theorem 4.3}
\newtheorem*{2glob}{Theorem 4.8}
\newtheorem*{3loc}{Theorem 5.3}
\newtheorem*{3glob}{Theorem 5.4}

\def\cal{{\rm cal}}
\def\cyl{{\rm cyl}}
\def\diam{{\rm diam}}
\def\Diff{{\rm Diff}}
\def\dist{{\rm dist}}
\def\Homeo{{\rm Homeo}}
\def\id{{\rm id}}
\def\Isom{{\rm Isom}}
\def\pr{{\rm pr}}
\def\st{{\rm st}}
\def\St{{\rm St}}

\numberwithin{equation}{section}

\begin{document}


\title[Sections of Serre fibrations with low-dimensional fibers]
{Sections of Serre fibrations with low-dimensional fibers}
\author{N. Brodsky}
\address{Department of Mathematics and Statistics,
University of Saskatche\-wan,
McLean Hall, 106 Wiggins Road, Saskatoon, SK, S7N 5E6, Canada}
\email{brodsky@math.usask.ca}
\author{A. Chigogidze}
\address{Department of Mathematics and Statistics,
University of Saskatche\-wan,
McLean Hall, 106 Wiggins Road, Saskatoon, SK, S7N 5E6, Canada}
\email{chigogid@math.usask.ca}
\author{E.V.~Shchepin}
\address{Steklov Institute of Mathematics,
Russian Academy of Science, Moscow 117966, Russia}
\email{scepin@mi.ras.ru}
\thanks{The second author was partially supported by NSERC research grant.}
\thanks{The third author was partially supported by Russian Foundation of Basic
Research (project 02-01-00014)}


\begin{abstract}{It was proved by H. Whitney in 1933 that it is possible to
mark a point in all curves in a continuous way. The main result of
this paper extends  the Whitney theorem to dimensions 2 and 3.
Namely, we prove that it is possible to choose a
point continuously in all two-dimensional surfaces sufficiently close to a
given surface, and in all $3$-manifolds sufficiently close to a
given 3-manifold.}
\end{abstract}

\maketitle
\markboth{N.~Brodsky, A.~Chigogidze, E.V.~Shchepin}
{Sections of Serre fibrations with low-dimensional fibers}

\section{Introduction}\label{S:intro}

It was proved by H. Whitney in 1933 \cite{Whitney} that it is possible to
mark a point in all curves in a continuous way. The main result of
this paper extends  the Whitney theorem to dimensions 2 and 3. To
be precise we prove that it is possible to choose a
point continuously in all two-dimensional surfaces sufficiently close to a
given surface, and in all $3$-manifolds sufficiently close to a
given 3-manifold. On the other hand this theorem fails for
dimensions greater than 4. In dimension 4 the question is an open problem.

The main tool used by Whitney is Whitney function. This
function produces a simultaneous parameterization of curves also known as the Morse parameterization. E.V.~Shchepin  conjectured that the analog of the Morse parameterization exists in dimensions less than 5. The
results of this paper count in favor of Shchepin's Conjecture in
dimensions 2 and 3.

Under the space $[M]$ of all manifolds of a given topological type
we mean the space of all closed subsets of a Hilbert cube which
are homeomorphic to the given manifold $M$ equipped with
Freshet topology. Let us denote by $[M,pt]$ the space of all
marked manifolds of the type $M$. Then there is a natural mapping
$p\colon [M,pt]\to[M]$ which is the universal Serre fibration with
fiber $M$. And  continuous choice of a point means continuous
section in this fibration. That is why we will speak later on
the sections of Serre fibrations only.

Now we can formulate

\smallskip

{\bf Shchepin's Conjecture.} \
{\it A Serre fibration with a locally arcwise connected metric base
is locally trivial if it has a low-dimensional $($of dimension $n \le 4)$
compact manifold as a constant fiber.}

\smallskip

In dimension $n=1$ the Shchepin's Conjecture is proved
even for non-compact fibers~\cite{RSS}.
Shchepin has proved \cite{S},\cite{DS} that positive solution of this
Conjecture in dimension $n$ implies positive solutions of both
CE-problem and Homeomorphism Group problem in dimension $n$. Since
$CE$-problem was solved in a negative way by A.N.~Dranishnikov, there
are dimensional restrictions in Shchepin's Conjecture.

In dimension $n=2$ we generalize the main result from~\cite{B}.

\begin{2loc}
Let $p\colon E\to B$ be a Serre fibration of $LC^0$-compacta
with a constant fiber which is a compact two-dimensional manifold.
If $B\in ANR$, then any section of $p$ over closed subset $A\subset B$
can be extended to a section of $p$ over some neighborhood of $A$.
\end{2loc}

In dimension $n=3$ we consider topologically regular mappings.
Note that if the Poincare Conjecture is true, then
any Serre fibration of $LC^2$-compacta
with a constant fiber which is a compact 3-manifold
is topologically regular~\cite{H}.

\begin{3loc}
Let $p\colon E\to B$ be a topologically regular mapping
of compacta with fibers homeomorphic to a 3-dimensional manifold.
If $B\in ANR$, then any section of $p$ over closed subset $A\subset B$
can be extended to a section of $p$ over some neighborhood of $A$.
\end{3loc}

Also we prove two theorems on global sections in Serre fibrations.

\begin{2glob}
Let $p\colon E\to B$ be a Serre fibration of $LC^0$-compactum $E$
onto an $ANR$-compactum $B$
with a constant fiber which is a connected two-dimensional compact manifold
$M$ not homeomorphic to the sphere or the projective plane.
Then $p$ admits a global section if either of the following conditions hold:
\begin{itemize}
\item[(a)] $\pi_1(M)$ is abelian and $H^2(B;\pi_1(F_b))=0$
\item[(b)] $\pi_1(M)$ is non-abelian, $M$ is not homeomorphic to the Klein bottle
   and $\pi_1(B)=0$.
\item[(c)] $M$ is homeomorphic to the Klein bottle and $\pi_1(B)=\pi_2(B)=0$.
\end{itemize}
\end{2glob}

\begin{3glob}
Let $p\colon E\to B$ be topologically regular mapping of compacta
with fibers homeomorphic to some compact connected 3-dimensional manifold $M$.
If $B$ is ANR-space, then $p$ admits a global section
if either of the following conditions hold:
\begin{itemize}
\item[(a)] $\pi_1(M)$ is abelian, $M$ is aspheric, and $H^2(B;\pi_1(F_b))=0$
\item[(b)] $M$ is closed hyperbolic 3-manifold and $\pi_1(B)=0$.
\item[(c)] $M$ is closed, irreducible, sufficiently large, contains
no embedded $\mathbb RP^2$ having a trivial normal bundle, and
$\pi_1(B)=\pi_2(B)=0$.
\end{itemize}
\end{3glob}

Let us recall some definitions and introduce our notations.
All spaces will be separable metrizable.
If not otherwise stated, by mapping we mean continuous single-valued mapping.
We equip the product $X\times Y$ with the metric
$$\dist_{X\times Y}((x,y),(x',y'))=\dist_X(x,x')+\dist_Y(y,y').$$
By $O(x,\varepsilon)$ we denote the open $\varepsilon$-neighborhood of the point $x$.

A multivalued mapping $F\colon X\to Y$ is called
{\it submapping} (or {\it selection}) of multivalued mapping
$G\colon X\to Y$ if $F(x)\subset G(x)$ for every $x\in X$.
The {\it gauge} of a multivalued mapping $F\colon X\to Y$ is defined as
$\cal (F)=\sup\{\diam F(x)\mid x\in X\}$.
The {\it graph} of multivalued mapping $F\colon X\to Y$
is the subset $\Gamma_F=\{(x,y)\in X\times Y\mid y\in F(x)\}$
of the product $X\times Y$.
For arbitrary subset $\mathcal U\subset X\times Y$ denote by
$\mathcal U(x)$ the subset $pr_Y(\mathcal U\cap(\{x\}\times Y))$ of $Y$.
Then for the graph $\Gamma_F$ we have $\Gamma_F(x)=F(x)$.

A multivalued mapping $G\colon X\to Y$ is called
{\it complete} if all sets $\{x\}\times G(x)$ are closed
with respect to some $G_\delta$-set $S\subset X\times Y$
containing the graph of this mapping.
A multivalued mapping $F\colon X\to Y$ is called
{\it upper semicontinuous} if for any open set
$U\subset Y$ the set $\{x\in X\mid F(x)\subset U\}$ is open in $X$.
A {\it compact} mapping is an upper semicontinuous
multivalued mapping with compact images of points.

An increasing\footnote{We   consider  only  increasing  filtrations
indexed by a segment of the natural series starting at zero.}
sequence (finite or infinite) of subspaces
$$ Z_0\subset Z_1\subset Z_2\subset\dots\subset Z $$
is called a {\it  filtration}  of space $Z$.
A sequence of multivalued mappings $\{F_k\colon X\to Y\}$
is called a {\it filtration of multivalued mapping} $F\colon X\to Y$
if for any $F(x)$, $\{F_k(x)\}$ is a filtration.

We say that a filtration of multivalued mappings $G_i\colon X\to Y$
is {\it complete} (resp. {\it compact}) if every mapping $G_i$
is complete (resp. compact).

\section{Local properties of multivalued mappings}\label{S:local}

Let $\gamma$ be a property of a topological space such that
every open subset inherits this property: if a space $X$
satisfies $\gamma$, then any open subset $U\subset X$
also satisfies $\gamma$.
We say that a space $Z$ satisfies $\gamma$ {\it locally}
if every point $z\in Z$ has a neighborhood with this property.

For a multivalued map $F\colon X\to Y$ to satisfy $\gamma$ locally
we not only require that every point-image $F(x)$ has this property locally,
but for any points $x\in X$ and $y\in F(x)$
there exist a neighborhood $W$ of $y$ in $Y$ and $U$ of $x$ in $X$
such that $W\cap F(x')$ satisfies $\gamma$ for every point $x'\in U$.
And we use the word "equi" for local properties of multivalued maps.

First example of such property is the local compactness.

\begin{defin}
A space $X$ is called {\it locally compact} if every point $x\in X$
has a compact neighborhood. We say that a multivalued map $F\colon X\to Y$ is {\it equi locally compact}
if for any points $x\in X$ and $y\in F(x)$
there exists a neighborhood $W$ of $y$ in $Y$ and $U$ of $x$ in $X$
such that $W\cap F(x')$ is compact for every point $x'\in U$.
\end{defin}

Another local property we are going to use is the hereditary asphericity.
This property is important in Geometric Topology (see \cite{DD})
and we will use the fact (which is easy to prove)
that the 2-dimensional Euclidean space is hereditary aspheric.
Recall that a compactum $K$ is called {\it approximately aspheric}
if for any (equivalently, for some) embedding of $K$ into an ANR-space $Y$
every neighborhood $U$ of $K$ in $Y$ contains a neighborhood $V$
with the following property: any mapping of the sphere $S^n$
into $V$ is homotopically trivial in $U$ provided $n\ge 2$.

\begin{defin}
We call a space $Z$ {\it hereditarily aspheric}
if any compactum $K\subset Z$ is approximately aspheric.

A space $Z$ is said to be {\it locally hereditarily aspheric}
if any point $z\in Z$ has a hereditarily aspheric neighborhood.
\end{defin}

Note that any 2-dimensional manifold is locally hereditarily aspheric.
Since 3-dimensional Euclidean space is not hereditarily aspheric, we
introduce a new property called hereditarily coconnected asphericity
to apply our technique to 3-dimensional manifolds.

\begin{defin}
We call a space $Z$ {\it hereditarily coconnectedly aspheric}
if any non-separating compactum $K\subset Z$ is approximately aspheric.

A space $Z$ is said to be {\it locally hereditarily coconnectedly aspheric}
if any point $z\in Z$ has a hereditarily coconnectedly aspheric neighborhood.
\end{defin}

Very important example of hereditarily coconnectedly aspheric space is
Euclidean 3-space~\cite{DS}. Therefore, any 3-dimensional manifold
is locally hereditarily coconnectedly aspheric.

\begin{defin}
We say that a multivalued map $F\colon X\to Y$ is
{\it equi locally hereditarily (coconnectedly) aspheric}
if for any points $x\in X$ and $y\in F(x)$ there exists a neighborhood $W$
of $y$ in $Y$ and $U$ of $x$ in $X$ such that $W\cap F(x')$ is
hereditarily (coconnectedly) aspheric for every point $x'\in U$.
\end{defin}

Now we consider different properties of pairs of spaces and define the
corresponding local properties for spaces and multivalued maps.
We follow definitions and notations from~\cite{DM}.

\begin{defin}
An ordering $\alpha$ of the subsets of a space $Y$ is {\it proper} provided:
\begin{itemize}
\item[(a)] If $W\alpha V$, then $W\subset V$;
\item[(b)] If $W\subset V$, and $V\alpha R$, then $W\alpha R$;
\item[(c)] If $W\alpha V$, and $V\subset R$, then $W\alpha R$.
\end{itemize}
\end{defin}

\begin{defin}
Let $\alpha$ be a proper ordering.
\begin{itemize}
\item[(a)]
A space $Y$ is {\it locally of type $\alpha$} if, whenever
$y\in Y$ and $V$ is a neighbourhood of $y$, then there a neighbourhood
$W$ of $y$ such that $W\alpha V$.
\item[(b)]
A multivalued mapping $F\colon X\to Y$ is {\it lower $\alpha$-continuous}
if for any points $x\in X$ and $y\in F(x)$ and for any neighbourhood $V$
of $y$ in $Y$ there exist neighbourhoods $W$ of $y$ in $Y$
and $U$ of $x$ in $X$ such that $(W\cap F(x'))\alpha(V\cap F(x'))$
provided $x'\in U$.
\end{itemize}
\end{defin}

For example, if $W\alpha V$ means that $W$ is contractible in
$V$, then locally of type $\alpha$ means locally contractible.
Another topological property which arises in this manner
is $LC^n$ (where $W\alpha V$ means that every continuous mapping
of the $n$-sphere into $W$ is homotopic to a constant mapping in $V$)
and the corresponding lower $\alpha$-continuity of multivalued map
is called {\it lower $(n+1)$-continuity}.
For the special case $n=-1$ the property $W\alpha V$ means that $V$ is
non-empty, and lower $\alpha$-continuity is the {\it lower semicontinuity}.

The following result is weaker than Lemma~3.5 from~\cite{BCK}.
We will use it with different properties $\alpha$ in
Sections~\ref{33} and~\ref{34}.

\begin{lem}\label{L32}
Let a lower $\alpha$-continuous mapping $\Phi\colon X\to Y$ of
compactum $X$ to a metric space $Y$ contains a compact submapping $F$.
Then for any $\varepsilon>0$ there exists a positive number $\delta$
such that for every point $(x,y)\in O(\Gamma_F,\delta)$ we have
$(O(y,\delta)\cap\Phi(x))\alpha (O(y,\varepsilon)\cap\Phi(x))$.
\end{lem}

In order to use results from~\cite{SB} we need a local property
called polyhedral $n$-connectedness.
A pair of spaces $V\subset U$ is called {\it polyhedrally $n$-connected}
if for any finite $n$-dimensional polyhedron $M$ and its
closed subpolyhedron $A$ any mapping of $A$ in $V$ can be extended
to a map of $M$ into $U$.
Note that for spaces being locally polyhedrally $n$-connected
is equivalent to be $LC^{n-1}$ (it follows from Lemma~\ref{lempolcont}).
The corresponding local property for multivalued map
is called {\it polyhedral lower $n$-continuity}.

\begin{lem} \label{lempolcont}
Any lower $n$-continuous multivalued mapping is
lower polyhedrally $n$-continuous.
\end{lem}

\begin{proof}
The proof easely follows from the fact that in connected filtration
$Z_0\subset Z_1\subset \dots\subset Z_n$ of spaces the pair $Z_0\subset Z_n$
is polyhedrally $n$-connected.
Given a mapping $f\colon A\to Z_0$ of subpolyhedron $A$ of $n$-dimensional
polyhedron $P$, we extend it successively over skeleta $P^{(k)}$
of $P$ such that the image of $k$-dimensional skeleton $P^{(k)}$
is contained in $Z_k$. Resulting map gives us an extension
$\widetilde f\colon P\to Z_n$ of $f$ which proves that the pair $Z_0\subset Z_n$
is polyhedrally $n$-connected.
\end{proof}

A filtration of multivalued maps $\{F_i\}$ is called
{\it polyhedrally connected} if every pair $F_{i-1}(x)\subset F_{i}(x)$
is polyhedrally $i$-connected.
A filtration $\{F_i\}$ is called {\it lower continuous}
if for any $i$ the mapping $F_i$ is lower $i$-continuous.

The following Lemma explains the reason to introduce the notion
of polyhedral $n$-connectivity.
This Lemma is a weak form of Compact Filtration Lemma from~\cite{SB}.

\begin{lem}\label{L31}
Any polyhedrally connected lower continuous finite filtration
of complete mappings of a compact space contains a compact approximately
connected\footnote{See Definition~\ref{appcon}.}
subfiltration of the same length.
\end{lem}

\begin{lem} \label{lemma2fibrprop}
If $p\colon E\to B$ is a Serre fibration of $LC^0$-compacta
with fibers homeomorphic to some compact 2-dimensional manifold, then
the multivalued mapping $p^{-1}\colon B\to E$ is
\begin{itemize}
\item equi locally hereditarily aspheric
\item polyhedrally lower $2$-continuous
\end{itemize}
\end{lem}

\begin{proof}
Since every open proper subset of a two-dimensional manifold is aspheric,
every compact proper subset of 2-manifold is approximately aspheric.
Therefore, the mapping $F$ is equi locally hereditarily aspheric.

It follows from a theorem of McAuley \cite{McAT} that the mapping $p^{-1}$ is
lower 2-continuous.
By Lemma~\ref{lempolcont}, the mapping $p^{-1}$ is
polyhedrally lower 2-continuous.
\end{proof}

We say that a subset $A$ of a space $Z$ is {\it coconnected}
if the complement $Z\setminus A$ is connected.
The following definition extends this property to pairs.

\begin{defin}
A pair $V\subset U$ of proper subsets of a space $Z$
is called {\it coconnected} if there exists a connected component of
$Z\setminus V$ containing $Z\setminus U$.
\end{defin}

If a pair $G_0\subset G_1$ of proper subsets of a space $Z$ is coconnected,
then we can define an operation of {\it coconnectification}
on subsets of $G_0$ as follows: for a subset $F_0\subset G_0$
its coconnectification is the union of $F_0$ and all components of
$Z\setminus F_0$ which do not intersect $Z\setminus G_1$.
Clearly, for a connected space $Z$ the coconnectification
of $F_0$ is the minimal subset
$F_1\subset G_1$ containing $F_0$ such that $Z\setminus F_1$ is connected.

\begin{defin}
A multivalued mapping $F\colon X\to Y$ is called {\it lower coconnected}
if for any points $x\in X$ and $y\in F(x)$ and for any sufficiently small neighbourhood $V$
of $y$ in $Y$ there exist neighbourhoods $W$ of $y$ in $Y$
and $U$ of $x$ in $X$ such that the pair $(W\cap F(x'))\subset(V\cap F(x'))$
is coconnected in $F(x')$ for every point $x'\in U$.
\end{defin}

If a multivalued mapping $F\colon X\to Y$ contains proper submappings
$G_0$ and $G_1$ such that for any $x\in X$ the pair $G_0(x)\subset G_1(x)$
is coconnected in $F(x)$, then for any submapping $F_0\subset G_0$
we define a {\it coconnectification} of $F_0$ as a multivalued mapping
taking a point $x\in X$ to the coconnectification of $F_0(x)$ in $G_1(x)$.

\begin{lem} \label{lemmacocon}
Suppose that lower $1$-continuous multivalued mapping $F\colon X\to Y$
contains proper submappings $G_0\subset G_1$ such that $G_1$ is compact and
for any $x\in X$ the pair $G_0(x)\subset G_1(x)$ is coconnected in $F(x)$.
Then for any compact submapping $F_0\subset G_0$ its coconnectification
$F_1$ is a compact submapping of $G_1$.
\end{lem}

\begin{proof}
Since the coconnectification of the set $F_0(x)$ is closed in $F(x)$
and is contained in $G_1(x)$, the set $F_1(x)$ is compact.

Let us prove that $F_1$ is u.s.c.
Suppose that for some point $y\in G_1(x)\setminus F_1(x)$ there is
a sequence $\{y_i\}_{i=1}^\infty$ of point converging to $y$ such that
$y_i$ belongs to a set $F_1(x_i)\setminus F_0(x_i)$ for some $x_i\in X$.
Fix a point $z\in F(x)\setminus G_1(x)$.
The points $y$ and $z$ belong to connected set $F(x)\setminus F_1(x)$
which is open in $F(x)$ and therefore is locally path connected
(since $F(x)\in LC^0$).
Hence, there exists a path $s\colon [0,1]\to F(x)\setminus F_1(x)$
such that $s(0)=y$ and $s(1)=z$.
Since $F$ is l.s.c. and $G_1$ is u.s.c., there is a sequence of points
$\{z_i\in F(x_i)\setminus G_1(x_i)\}_{i=M}^\infty$, converging to $z$.

Using lower $1$-continuity of the mapping $F$
we can choose a sequence of maps
$\{s_i\colon [0,1]\to F(x_i)\}_{i=M'}^\infty$ such that
$s_i(0)=y_i$, $s_i(1)=z_i$ and $s_i\to_{i\to\infty} s$.
Since the path $s$ does not intersect the fiber $F_0(x)$
and $F_0$ is u.s.c., all but the finite number of paths $s_i$
do not intersect the fibers of $F_0$.
It means that the points $y_i$ and $z_i$ belong to
the same connected component of the set $F(x_i)\setminus F_0(x_i)$,
which contradicts to the choices of $y_i$ and $z_i$.
\end{proof}

\begin{defin}
The mapping $f \colon X \to Y$ is said to be {\it topologically regular}
provided that if $\varepsilon > 0$ and $y \in Y$,
then there is a positive number $\delta$ such that
$dist(y,y') < \delta$,  $y' \in Y$,  implies that
there is a homeomorphism of $f^{-1}(y)$ onto $f^{-1}(y')$
which moves no point as much as $\varepsilon$
(i.e. an $\varepsilon$-homeomorphism).
\end{defin}

Note that if the Poincare Conjecture is true, then
any Serre fibration $f \colon X \to Y$ of $LC^2$-compacta
with a constant fiber which is a compact three-dimensional manifold
is topologically regular~\cite{H}.

\begin{lem} \label{lemmafibrprop}
If $p\colon E\to B$ is a topologically regular mapping of compacta
with fibers homeomorphic to some connected 3-dimensional manifold, then
the multivalued mapping $p^{-1}\colon B\to E$ is
\begin{itemize}
\item equi locally hereditarily coconnectedly aspheric
\item lower coconnected
\item equi locally compact
\item lower polyhedrally $2$-continuous
\end{itemize}
\end{lem}

\begin{proof}
Clearly, we can identify the graph of $p^{-1}$ with the space $E$
and the projection of $\Gamma_{p^{-1}}$ onto $B$ with $p$.
Fix a point $q\in E$ and $\varepsilon>0$.
We will find $\delta>0$ such that for any point $x\in p(O(q,\delta))$
there exist subsets $D^3$ and $O^3$ of the fiber $p^{-1}(x)$ such that
\[  O(q,\delta)\cap p^{-1}(x)\subset D^3\subset O^3\subset
    O(q,\varepsilon)\cap p^{-1}(x)  \]
where $D^3$ is homeomorphic to closed 3-ball and
$O^3$ is homeomorphic to $\mathbb R^3$.
Then last three properties of $p^{-1}$ follow easely.
And the first property follows from Lemma~2.4 from~\cite{DS}
which states that a compactum in $\mathbb R^3$ which does not separate
$\mathbb R^3$ is approximately aspheric.

Take a neighborhood $O^3_q$ of the point $q$ in the fiber $p^{-1}(p(q))$ such
that $O^3_q$ is homeomorphic to $\mathbb R^3$ and is contained in $O(q,\varepsilon/2)$.
Note that if $h$ is $\varepsilon/2$-homeomorphism of $O^3_q$, then
$h(O^3_q)$ is contained in $O(q,\varepsilon)$.
Let $D^3_q$ be a neighborhood of $q$ in $p^{-1}(p(q))$
homeomorphic to closed 3-ball.
Take a number $\sigma>0$ such that $O(q,\sigma)\cap p^{-1}(p(q))$
is contained in $D^3_q$.
Choose a positive number $\delta<\sigma/2$ such that for any point
$x\in O(p(q),\delta)$ there exists $\sigma/2$-homeomorphism
of the fiber $p^{-1}(p(q))$ onto $p^{-1}(x)$.
Now take a point $x\in p(O(q,\delta))$ and fix $\sigma/2$-homeomorphism $h$
of the fiber $p^{-1}(p(q))$ onto $p^{-1}(x)$.
By the choice of $\sigma$, the set $h(D^3_q)$ contains
$O(q,\sigma/2)\cap p^{-1}(x)$. Therefore, we have
\[  O(q,\delta)\cap p^{-1}(x)\subset h(D^3_q)\subset h(O^3_q)\subset
    O(q,\varepsilon)\cap p^{-1}(x).  \]
\end{proof}

\section{Singlevalued approximations}\label{S:approx}

In this Section we prove filtered finite dimensional approximation theorem
(Theorem~\ref{thmappr}) and then apply it in a usual way (compare~\cite{GGK})
to prove an approximation theorem for maps of ANR-spaces.
Since we are going to use singular filtrations of multivalued maps
instead of usual filtrations, our Theorem~\ref{thmappr} generalizes
the Filtered approximation theorem proved in~\cite{SB}.
But the proof of our singular version of
filtered approximation theorem in full generality requires
a lot of technical details to establish.
So, we decided to prove precisely the version that we need ---
for compact maps of metric spaces.

Let us introduce a notion of singular filtration.

\begin{defin}
A {\it singular pair} of spaces is a triple $(Z,\phi,Z')$
where $\phi\colon Z\to Z'$ is a mapping.

We say that a space $Z$ contains a {\it singular filtration
of spaces} if a finite sequence of pairs $\{(Z_i,\phi_i)\}_{i=0}^n$
is given where $Z_i$ is a space and $\phi_i\colon Z_i\to Z_{i+1}$ is a map
(we identify $Z_{n+1}$ with $Z$).
\end{defin}

For a multivalued map $F\colon X\to Y$ it is useful to consider its
{\it graph fibers} $\{x\}\times F(x)\subset \Gamma_F$
instead of usual fibers $F(x)\subset Y$.
While the graph fibers are always homeomorphic to the usual fibers,
different graph fibers do not intersect (the usual fibers
may intersect in $Y$). We denote the graph fiber of the map $F$
over a point $x\in X$ by $F^\Gamma(x)$.

To define the notion of singular filtration for multivalued maps
we introduce a notion of fiberwise transformation of multivalued maps.

\begin{defin}
For multivalued mappings $F$ and $G$ of a space $X$ a
{\it fiberwise transformation from $F$ to $G$}
is a continuous mapping $T\colon \Gamma_F\to \Gamma_G$
such that $T(F^\Gamma(x))\subset G^\Gamma(x)$ for every $x\in X$.

A {\it fiber} $T(x)$ of the fiberwise transformation $T$ over
the point $x\in X$ is a mapping $T(x)\colon F(x)\to G(x)$
determined by $T$.

We say that a multivalued mapping $F\colon X\to Y$ contains a
{\it singular filtration of multivalued maps}
if a finite sequence of pairs $\{(F_i,T_i)\}_{i=0}^n$
is given where $F_i\colon X\to Y_i$ is a multivalued mapping
and $T_i$ is a fiberwise transformation from $F_i$ to $F_{i+1}$
(we identify $F_{n+1}$ with $F$).
\end{defin}

To construct continuous approximations of multivalued maps we
need the notion of approximate asphericity.

\begin{defin}
A pair of compacta $K\subset K'$ is called
{\it approximately $n$-aspheric} if for any embedding of
$K'$ into $ANR$-space $Z$
for every neighborhood $U$ of $K'$ in $Z$ there exists
a neighborhood $V$ of $K$ such that any mapping
$f\colon S^n\to V$ is homotopically trivial in $U$.

A compactum $K$ is {\it approximately $n$-aspheric} if the pair
$K\subset K$ is approximately $n$-aspheric.
\end{defin}

The following is a singular version of approximate asphericity.

\begin{defin}
A singular pair of compacta $(K,\phi,K')$ is called
{\it approximately $n$-aspheric} if for any embeddings
$K\subset Z$ and $K'\subset Z'$ in $ANR$-spaces and for
any extension of $\phi$ to a map $\widetilde \phi\colon OK\to Z'$
of some neighborhood $OK$ of $K$ the following holds:
for every neighborhood $U$ of $K'$ in $Z'$ there exists
a neighborhood $V$ of $K$ in $OK$ such that for any mapping
$f\colon S^n\to V$ the spheroid $\widetilde\phi\circ f\colon S^n\to U$
is homotopically trivial in $U$.
\end{defin}

Following R.C.~Lacher~\cite{L}, one can prove that this notion does
not depend on the choices of $ANR$-spaces $Z$ and $Z'$
and on the embeddings of $K$ and $K'$ into these spaces.

\begin{defin}
A singular filtration of compacta $\{(K_i,\phi_i)\}_{i=0}^n$
is called {\it approximately connected} if for every $i<n$
the singular pair $(K_i,\phi_i,K_{i+1})$ is approximately $i$-aspheric.
\end{defin}

Clearly, a singular pair of compacta $(K,\phi,K')$ is approximately
$n$-aspheric in either of the following three situations:
compactum $K$, compactum $K'$, or the pair $\phi(K)\subset K'$
is approximately $n$-aspheric.

\begin{defin}\label{appcon}
A singular filtration $\mathcal F=\{(F_i,T_i)\}_{i=0}^n$ of compact
mappings $F_i\colon X\to Y_i$ is said to be {\it approximately connected}
if for every point $x\in X$ the singular filtration of compacta
$\{(F_i(x),T_i(x))\}_{i=0}^n$ is approximately connected.

An approximately connected singular filtration
$\mathcal F=\{(F_i\colon X\to Y_i,T_i)\}_{i=0}^n$
is said to be {\it approximately $\infty$-connected}
if the mapping $F_n$ has approximately $k$-aspheric point-images
$F_n(x)$ for all $k\ge n$ and all $x\in X$.
\end{defin}

Note that if a singular filtration $\mathcal F=\{(F_i,T_i)\}_{i=0}^n$
is approximately $\infty$-connected, then the mapping $F_n$ contains
an approximately connected singular filtration of any given finite length.

We will reduce our study of singular filtrations to the study
of usual filtrations using the following cylinder construction.

\begin{defin}
For a continuous singlevalued mapping $f\colon X\to Y$ we define a
{\it cylinder} of $f$ denoted by $\cyl (f)$ as a space obtained
from the disjoint union of $X\times [0,1]$ and $Y$ by
identifying each $\{x\}\times \{1\}$ with $f(x)$.
\end{defin}

Note that the cylinder $\cyl (f)$ contains a homeomorphic copy of $Y$
called the {\it bottom} of the cylinder, and a homeomorphic copy of $X$
as $X\times \{0\}$ called the {\it top} of the cylinder.

\begin{rem} \label{remretr}
There is a natural deformation retraction $r\colon \cyl(f)\to Y$
onto the bottom $Y$. Clearly, the fiber of the mapping $r$
over a point $y\in Y$ is either one point $\{y\}$
or a cone over the set $f^{-1}(y)$.
Therefore, if the map $f$ is proper, then $r$ is $UV^\infty$-mapping.
\end{rem}

\begin{rem} \label{rememb}
Suppose that $X$ is embedded into Banach space $B_1$ and
$Y$ is embedded into Banach space $B_2$.
Then we can naturally embed the cylinder $\cyl (f)$ into the
product $B_1\times \Bbb R\times B_2$. The embedding is clearly defined
on the top as embedding into $B_1\times \{0\}\times\{0\}$
and on the bottom as embedding into $\{0\}\times \{1\}\times B_2$.
We extend these embeddings to the whole cylinder by sending its point
$\{x\}\times \{t\}$ to the point
$\{(1-t)\cdot x\}\times t\times \{t\cdot f(x)\}$.
\end{rem}

\begin{lem} \label{lemasphcomp}
If a singular pair of compacta $(K,\phi,K')$ is approximately $n$-aspheric,
then the pair $K\subset \cyl(\phi)$ is approximately $n$-aspheric.
\end{lem}

\begin{proof}
Let us fix embeddings of $K$ into Banach space $B_1$,
of $K'$ into Banach space $B_2$, and of the cylinder $\cyl(\phi)$
into the product $B=B_1\times \Bbb R\times B_2$ as described
in Remark~\ref{rememb}.
Fix a neighborhood $U$ of $\cyl(\phi)$ in $B$.
Extend the mapping $\phi$ to a map $\phi_1\colon B_1\to B_2$.
Take a neighborhood $V_1$ of the top of our cylinder in $B_1$
such that the cylinder $\cyl(\phi_1|_{V_1})$ is contained in $U$.
Using approximate $n$-asphericity of the pair $(K,\phi,K')$
we find for a neighborhood $U\cap \{0\}\times\{1\}\times B_2$ of $K'$ in
$\{0\}\times\{1\}\times B_2$ a neighborhood $V'$ of $K$ in
$B_1\times\{0\}\times\{0\}$.
Let $\varepsilon$ be a positive number such that the product
$V=V'\times (-\varepsilon,\varepsilon)\times O(0,\varepsilon)$ is contained in $U$.

Given a spheroid $f\colon S^n\to V$ we retract it into
$V'\times\{0\}\times\{0\}$, then retract it to the bottom
of the cylinder $\cyl(\phi_1|_{V_1})$ using Remark~\ref{remretr},
and finally contract it to a point inside $U\cap \{0\}\times\{1\}\times B_2$.
Clearly, the whole retraction sits inside $U$, as required.
\end{proof}

\begin{defin}
Let $\mathcal F=\{(F_i\colon X\to Y_i,T_i)\}_{i=0}^n$
be a singular filtration of a multivalued mapping $F\colon X\to Y=Y_{n+1}$.
If all the spaces $Y_i$ are Banach, then for a multivalued mapping $\mathbb F$
from $X$ to $\mathbb Y=Y\times \prod_{i=0}^n (Y_i\times\mathbb R)$
defined as $\mathbb F(x)=\cup_{k=0}^n\cyl(T_k(x))$
we can define a {\it cylinder} $\cyl(\mathcal F)$ as a filtration
of multivalued maps $\{\mathbb F_i\}_{i=0}^n$ defined as follows:
\[  \mathbb F_0=F_0 \text{ \ and \ }
    \mathbb F_i(x)=\bigcup_{k=0}^{i-1}\cyl(T_k(x)).  \]
\end{defin}

It is easy to see that for a singular filtration
$\mathcal F=\{(F_i,T_i)\}_{i=0}^n$ of compact mappings $F_i$
the filtration $\cyl(\mathcal F)$ consists of compact mappings $\mathcal F_i$.

\begin{lem} \label{lemapprfil}
If a singular filtration $\mathcal F=\{(F_i,T_i)\}_{i=0}^n$ of compact maps
is approximately connected, then the filtration
$\cyl(\mathcal F)=\{\mathbb F_i\}_{i=0}^n$ is approximately connected.
\end{lem}

\begin{proof}
Using Remark~\ref{remretr} it is easy to define a deformation retraction
$r\colon \mathbb F_i(x)\to F_i(x)$ which is $UV^\infty$-mapping.
This retraction defines $UV^\infty$-mapping of pairs
$(\mathbb F_i(x),\mathbb F_{i+1}(x))\to (F_i(x),\cyl(T_i(x)))$.
By Lemma~\ref{lemasphcomp} the pair $F_i(x)\subset\cyl(T_i(x))$
is approximately $i$-aspheric, and by Pairs Mapping Lemma from~\cite{SB}
the pair $\mathbb F_i(x)\subset\mathbb F_{i+1}(x)$
is also approximately $i$-aspheric
\end{proof}

\begin{thm} \label{thmappr}
Let $H\colon X\to Y$ be a multivalued mapping of metric space $X$
to a Banach space $Y$.
If $\dim X\le n$ and $H$ contains approximately connected
singular filtration $\mathcal H=\{(H_i\colon X\to Y_i,T_i)\}_{i=0}^n$
of compact mappings, then any neighborhood $\mathcal U$
of the graph $\Gamma_H$ contains the graph of a single-valued
and continuous mapping $h\colon X\to Y$.
\end{thm}

\begin{proof}
Without loss of generality we may assume that all spaces $Y_i$
are Banach spaces. We consider $Y$ as a subspace of the product
$\mathbb Y=Y\times \prod_{i=0}^n (Y_i\times\mathbb R)$.
Clearly, $H$ is a submapping of a multivalued mapping
$\mathbb H\colon X\to \mathbb Y$ defined as
$\mathbb H(x)=\cup_{k=0}^n\cyl(T_k(x))$ and
$\Gamma_{\mathbb H}$ admits a deformation retraction $R$ onto $\Gamma_H$.
Fix a neighborhood $\mathcal U$ of the graph $\Gamma_H$ in $X\times Y$.
Since all maps $H_i$ are compact, $\mathbb H$ is also compact
and the graph $\Gamma_{\mathbb H}$ is closed in $X\times\mathbb Y$.
Extend the mapping $\pr_Y\circ R\colon \Gamma_{\mathbb H}\to Y$
to some neighborhood $\mathcal W$ of $\Gamma_{\mathbb H}$
in $X\times\mathbb Y$ and denote by $R'$ the map of $\mathcal W$
to $X\times Y$ such that $\pr_Y\circ R'$ is our extension.
Clearly, we may assume that $R'(\mathcal W)$ is contained in $\mathcal U$.

By Lemma~\ref{lemapprfil} the multivalued map $\mathbb H$
admits approximately connected filtration $\cyl(\mathcal H)$
of compact multivalued maps.
By Single-Valued Approximation Theorem from~\cite{SB} there exists
a singlevalued continuous mapping ${\bf h}\colon X\to \mathbb Y$
with $\Gamma_{\bf h}\subset\mathcal W$.
Define a singlevalued continuous map $h$ by the equality
$\Gamma_h=R'(\Gamma_{\bf h})$.
Clearly, $\Gamma_h$ is contained in $R'(\mathcal W)\subset \mathcal U$.
\end{proof}

\begin{thm} \label{thmapprox}
Suppose that a compact mapping of separable metric ANRs $F\colon X\to Y$
admits a compact singular approximately $\infty$-connected filtration.
Then for any compact space $K\subset X$ every neighborhood of
the graph $\Gamma_F(K)$ contains the graph of a single-valued
and continuous mapping $f\colon K\to Y$.
\end{thm}

\begin{proof}
Let $\mathcal U$ be an open neighborhood of the graph $\Gamma_F(K)$
in the product $X\times Y$. Since $F$ is upper semicontinuous,
there is a neighborhood $OK$ of compactum $K$ such that
$\Gamma_F(OK)$ is contained in $\mathcal U$.
Since any open subset of separable ANR-space is separable ANR-space
\cite{Han}, we can denote $OK$ by $X$ and consider $\mathcal U$
as an open neighborhood of the graph $\Gamma_F$.

For every point $x\in X$ take open neighborhoods $O_x\subset X$
of the point $x$ and $V_x\subset X$ of the compactum $F(x)$
such that the product $O_x\times V_x$ is contained in $\mathcal U$.
Using upper semicontinuity of $F$ we can choose $O_x$
so small that the following inclusion holds: $F(O_x)\subset V_x$.
Fix an open covering $\omega_1$ of the space $X$ which is
starlike refinement of $\{O_x\}_{x\in X}$.
Let $\omega_2$ be a locally finite open covering of the space $X$
which is starlike refined into $\omega_1$.

There exist a locally finite simplicial complex $L$
and mappings $r\colon X\to L$ and $j\colon L\to X$ such that
the map $j\circ r$ is $\omega_2$-close to $\id_X$ \cite{Han}.
Fix a finite subcomplex $N\subset L$ containing the compact set $r(K)$.
Define a compact mapping $\Psi\colon N\to Y$ by the formula $\Psi=F\circ j$.
Clearly, the mapping $\Psi$ admits a compact approximately connected
singular filtration of any length (particularly, of the length $\dim N$).
Let us define a neighborhood $\mathcal W$ of the graph $\Gamma_\Psi$.
For every point $q\in N$ we put

$$ \mathcal W(q)=\bigcap\{\mathcal U(y)\mid y\in\st_{\omega_1}(\St_{\omega_2}(j(q)))\}. $$

By $\St_{\omega_2}(j(q))$, we denote the star of the point $j(q)$
with respect to the covering $\omega_2$.
And by $\st(A,\omega)$, we denote the set
$\bigcup\{U\in\omega\mid A\subset U\}$.

By Theorem~\ref{thmappr} there exists a single-valued continuous mapping
$\psi\colon N\to Y$ such that the graph $\Gamma_\psi$ is contained in $\mathcal W$.
Put $f=\psi\circ r\colon X\to Y$.
For any point $x\in K$ we have
$\psi(r(x))\in\cap\{\mathcal U(x')\mid x'\in\St_{\omega_2}(j\circ r(x))\}$.
Since $x\in\St_{\omega_2}(j\circ r(x))$, then $\psi(r(x))\in\mathcal U(x)$.
That is, the graph of $f$ is contained in $\mathcal U$.
\end{proof}

\section{Fibrations with 2-manifold fibers}\label{33}

Our strategy of constructing a section of a Serre fibration is as follows.
We consider the inverse (multivalued) mapping and find
its compact submapping admitting continuous approximations.
Then we take very close continuous approximation and use it to
find again a compact submapping with small diameters of fibers
admitting continuous approximations.
When we continue this process we get a sequence of compact submappings
with diameters of fibers tending to zero.
This sequence will converge to the desired singlevalued submapping (selection).

\begin{lem} \label{lemmashrinking}
Let $F\colon X\to Y$ be equi locally hereditarily aspheric,
lower 2-continuous complete multivalued mapping
of $ANR$-space $X$ to Banach space $Y$.
Suppose that a compact submapping $\Psi\colon A\to Y$ of $F|_A$
is defined on a compactum $A\subset X$ and admits continuous approximations.
Then for any $\varepsilon>0$ there exists a neighborhood $OA$ of $A$
and a compact submapping $\Psi'\colon OA\to Y$ of $F|_{OA}$
such that $\Gamma_{\Psi'}\subset O(\Gamma_\Psi,\varepsilon)$,
$\Psi'$ admits a compact approximately $\infty$-connected
filtration, and $\cal \Psi'<\varepsilon$.
\end{lem}

\begin{proof}
Fix a positive number $\varepsilon$.
Apply Lemma~\ref{L32} with $\alpha$ being equi local hereditary asphericity
to get a positive number $\varepsilon_2<\varepsilon/4$.
By Lemma~\ref{lempolcont} the mapping $F$ is lower polyhedrally 2-continuous.
Subsequently applying Lemma~\ref{L32} with $\alpha$ being
polyhedral $n$-continuity for $n=2,1,0$, we find positive numbers
$\varepsilon_1$, $\varepsilon_0$, and $\delta$ such that $\delta<\varepsilon_0<\varepsilon_1<\varepsilon_2$
and for every point $(x,y)\in O(\Gamma_{\Psi},\delta)$
the pair $(O(y,\varepsilon_1)\cap F(x),O(y,\varepsilon_2)\cap F(x))$
is polyhedrally 2-connected,
the pair $(O(y,\varepsilon_0)\cap F(x),O(y,\varepsilon_1)\cap F(x))$
is polyhedrally 1-connected, and
the intersection $O(y,\varepsilon_0)\cap F(x)$ is not empty.

Let $f\colon K\to Y$ be a continuous single-valued mapping
whose graph is contained in $O(\Gamma_{\Psi},\delta)$.
Let $f'\colon\mathcal OK\to Y$ be a continuous extension of the mapping $f$
over some neighborhood $\mathcal OK$ such that the graph of $f'$
is contained in $O(\Gamma_{\Psi},\delta)$.
Now we can define a polyhedrally connected filtration
$G_0\subset G_1\subset G_2\colon \mathcal OK\to Y$
of the mapping $F|_{\mathcal OK}$ by the equality
\[
G_i(x)=O(f'(x),\varepsilon_i)\cap F(x).
\]
Since the set $\cup_{x\in\mathcal OK}\{\{x\}\times O(f'(x),\varepsilon_i)\}$
is open in the product $\mathcal OK\times Y$ and the mapping $F$
is complete, then $G_i$ is also complete.
Clearly, $\cal G_2<2\varepsilon_2<\varepsilon$ and for any point $x\in K$
the set $G_2^\Gamma(x)$ is contained in $O(\Gamma_\Psi,\varepsilon)$.
Now, applying Lemma~\ref{L31} to the filtration $G_0\subset G_1\subset G_2$,
we obtain a compact approximately connected subfiltration
$F_0\subset F_1\subset F_2\colon \mathcal OK\to Y$.
By the choice of $\varepsilon_2$ the mapping $F_2$ has approximately aspheric
point-images. Therefore, the filtration $F_0\subset F_1\subset F_2$
is approximately $\infty$-connected.
Finally, we put $\Psi'=F_2$.
\end{proof}

\begin{thm} \label{thmapproxsections}
Let $F\colon X\to Y$ be equi locally hereditarily aspheric,
lower 2-continuous complete multivalued mapping
of locally compact $ANR$-space $X$ to Banach space $Y$.
Suppose that a compact submapping $\Psi\colon A\to Y$ of $F|_A$
is defined on compactum $A\subset X$ and admits continuous approximations.
Then for any $\varepsilon>0$ there exists a neighborhood $OA$ of $A$
and a single-valued continuous selection $s\colon OA\to Y$
of $F|_{OA}$ such that $\Gamma_s\subset O(\Gamma_\Psi,\varepsilon)$.
\end{thm}

\begin{proof}
Consider a $G_\delta$-subset $G\subset X\times Y$ such that
all fibers of $F$ are closed in $G$ and fix open sets
$G_i\subset X\times Y$ such that $G=\cap^\infty_{i=1} G_i$.
Fix $\varepsilon>0$ such that $O(\Gamma_\Psi,\varepsilon)\subset G_1$.
By Lemma~\ref{lemmashrinking} there is a neighborhood $U_1$
of $A$ in $X$ and a compact submapping $\Psi_1\colon U_1\to Y$ of $F|_{U_1}$
such that $\Gamma_{\Psi_1}\subset O(\Gamma_\Psi,\varepsilon)$,
$\Psi_1$ admits a compact approximately $\infty$-connected
filtration, and $\cal \Psi_1<\varepsilon$.
Since $X$ is locally compact and $A$ is compact, there exists a compact
neighborhood $OA$ of $A$ such that $OA\subset U_1$.
By Theorem~\ref{thmapprox} the mapping $\Psi_1|_{OA}$
admits continuous approximations.
Take $\varepsilon_1<\varepsilon$ such that the neighborhood
$\mathcal U_1=O(\Gamma_{\Psi_1}(OA),\varepsilon_1)$ lies in $O(\Gamma_\Psi,\varepsilon)$.
Clearly, $\mathcal U_1\subset G_1$.

Now by induction with the use of Lemma~\ref{lemmashrinking},
we construct a sequence of neighborhoods
$U_1\supset U_2\supset U_3\supset\dots$ of the
compactum $OA$, a sequence of compact submappings
$\{\Psi_k\colon U_k\to Y\}^\infty_{k=1}$ of the mapping $F$,
and a sequence of neighborhoods
$\mathcal U_k=O(\Gamma_{\Psi_1}(OA),\varepsilon_k)$
such that for every $k\ge 2$ we have $\cal \Psi_k<\varepsilon_{k-1}/2<\varepsilon/2^k$,
and $\mathcal U_k(OA)$ is contained in $\mathcal U_{k-1}(OA)\cap G_k$.
It is not difficult to choose the neighborhood $\mathcal U_k$
of the graph $\Gamma_{\Psi_k}$ in such a way that for every point $x\in U_k$
the set $\mathcal U_k(x)$ has diameter less than $\frac{3}{2^k}$.

Then for any $m\ge k\ge 1$ and for any point $x\in OA$ we have
$\Psi_m(x)\subset O(\Psi_k(x),\frac{3}{2^k})$; this implies the fact
that the sequence $\{\Psi_k|_{OA}\}_{k=1}^\infty$ is a Cauchy sequence.
Since $Y$ is complete, there exists the limit $s\colon OA\to Y$ of this sequence.
The mapping $s$ is single-valued by
the condition $\cal\Psi_k<\frac{1}{2^k}$ and is upper
semicontinuous (and, therefore, is continuous) by the upper
semicontinuity of all the mappings $\Psi_k$.
Clearly, for any $x\in OA$ the point $s(x)$ lies in $G(x)$
and is a limit point of the set $F(x)$.
Since $F(x)$ is closed in $G(x)$, then $s(x)\in F(x)$,
i.e. $s$ is a selection of the mapping $F$.
\end{proof}

\begin{thm}
Let $p\colon E\to B$ be a Serre fibration of $LC^0$-compacta
with a constant fiber which is a compact two-dimensional manifold.
If $B\in ANR$, then any section of $p$ over closed subset $A\subset B$
can be extended to a section of $p$ over some neighborhood of $A$.
\end{thm}

\begin{proof}
Let $s\colon A\to E$ be a section of $p$ over $A$.
Embed $E$ into Hilbert space $l_2$ and consider a multivalued mapping
$F\colon B\to l_2$ defined as follows:
$$
   F(b)=\cases s(b), &\text{if \;$b\in A$}\\
               p^{-1}(b),   &\text{if \;$x\in B\setminus A$.}\endcases
$$
Since every fiber $p^{-1}(b)$ is compact, the mapping $F$ is complete.
By Lemma~\ref{lemma2fibrprop} the mapping $F$ is equi
locally hereditarily aspheric and lower $2$-continuous.
We can apply Theorem~\ref{thmapproxsections} to the mapping $F$
and its submapping $s$ to find a single-valued continuous selection
$\widetilde s\colon OA\to l_2$ of $F|_{OA}$.
By definition of $F$, we have $\widetilde s|_{A}=F|_{A}=s$.
Clearly, $\widetilde s$ defines a section of the fibration $p$
over $OA$ extending $s$.
\end{proof}

\begin{defin}
For a mapping $p\colon E\to B$ we say that $s\colon B\to E$ is
{\it $\varepsilon$-section} if the map $p\circ s$ is $\varepsilon$-close to
the identity $\id_B$.
\end{defin}

The following proposition easily follows from Theorem~4.1 of the paper~\cite{Mi}.

\begin{pro} \label{profibr}
If $p\colon E\to B$ is a locally trivial fibration of finite-dimensional
compacta with locally contractible fiber, then there is $\varepsilon>0$ such that an
existence of $\varepsilon$-section for $p$ implies an existence of a section for $p$.
\end{pro}

We will use the following two statements in the proof of
existence of global sections in Serre fibrations.
The proof of these two Propositions follows from the Bestvina's
construction of Menger manifold~\cite{bestvina}
and Dranishnikov's triangulation theorem for Menger manifolds~\cite{Dr2}.

\begin{pro} \label{propoly}
Let $X$ be a compact 2-dimensional Menger manifold.
For any $\varepsilon>0$ there exist a finite polyhedron $P$ and maps $g\colon X\to P$
and $h\colon P\to X$ such that $h\circ g$ is $\varepsilon$-close to the identity.
If $\pi_1(X)=0$, then we may choose $P$ with $\pi_1(P)=0$.
\end{pro}

\begin{pro} \label{pro3poly}
Let $X$ be a compact 3-dimensional Menger manifold with $\pi_1(X)=\pi_2(X)=0$.
For any $\varepsilon>0$ there exist a finite polyhedron $P$ with
$\pi_1(P)=\pi_2(P)=0$ and maps $g\colon X\to P$
and $h\colon P\to X$ such that $h\circ g$ is $\varepsilon$-close to the identity.
\end{pro}

\begin{thm} \label{glob2sec}
Let $p\colon E\to B$ be a Serre fibration of $LC^0$-compactum $E$
onto $ANR$-compactum $B$
with a constant fiber which is a compact connected two-dimensional manifold
$M$ not homeomorphic to sphere or projective plane.
Then $p$ admits a global section if either of the following conditions hold:
\begin{itemize}
\item[(a)] $\pi_1(M)$ is abelian and $H^2(B;\pi_1(F_b))=0$
\item[(b)] $\pi_1(M)$ is non-abelian, $M$ is not homeomorphic to the Klein bottle
   and $\pi_1(B)=0$.
\item[(c)] $M$ is homeomorphic to the Klein bottle and $\pi_1(B)=\pi_2(B)=0$.
\end{itemize}
\end{thm}

\begin{proof}
Embed $E$ into the Hilbert space $l_2$ and consider a multivalued mapping
$F\colon B\to l_2$ defined as $F=p^{-1}$.
Since every fiber $p^{-1}(b)$ is compact, the mapping $F$ is complete.
It follows from Lemma~\ref{lemma2fibrprop} that the mapping $F$ is
equi locally hereditarily aspheric and lower $2$-continuous.

Now we show that $F$ admits a compact singular approximately
$\infty$-connected filtration.
In cases (a) and (b) there exists $UV^1$-mapping $\mu$
of Menger 2-dimensional manifold $L$
onto $B$~\cite{Dr1}. Note that $\pi_1(L)=0$ if $\pi_1(B)=0$.
In case (c) we consider $UV^2$-mapping $\mu$
of Menger 3-dimensional manifold $L$ onto $B$~\cite{Dr1};
note that $\pi_1(L)=\pi_2(L)=0$ if $\pi_1(B)=\pi_2(B)=0$.
Since $\dim L<\infty$, the induced fibration $p_L=\mu^*(p)\colon E_L\to L$
is locally trivial~\cite{HD}.
By Proposition~\ref{profibr} there is $\varepsilon>0$ such that an
existence of $\varepsilon$-section for $p_L$ implies an
existence of a section for $p_L$.
In cases (a) and (b), by Proposition~\ref{propoly}, there exist a 2-dimensional
finite polyhedron $P$ and continuous maps $g\colon L\to P$
and $h\colon P\to L$ such that $h\circ g$ is $\varepsilon$-close
to the identity (we assume $\pi_1(P)=0$ in case $\pi_1(B)=0$).
In case (c) by Proposition~\ref{pro3poly} there exist a 3-dimensional
finite polyhedron $P$ with $\pi_1(P)=\pi_2(P)=0$
and continuous maps $g\colon L\to P$ and $h\colon P\to L$
such that $h\circ g$ is $\varepsilon$-close to the identity.

Consider a locally trivial fibration $p_P=h^*(p_L)\colon E_P\to P$.

{\bf Claim.} \ {\it The fibration $p_P$ has a section $s_P$.}

\begin{proof}
\begin{itemize}
\item[(a)]
If $\pi_1(M)$ is abelian and $H^2(B;\pi_1(F_b))=0=H^2(P;\pi_1(F_b))$,
the fibration $p_P$ has a section $s_P$~\cite{W}.

\item[(b)]
Since $\pi_1(P)=0$ and $\dim P=2$, then $P$ is homotopy equivalent
to a bouquet of 2-spheres $\Omega=\vee_{i=1}^m S^2_i$.
Let $\psi\colon P\to\Omega$ and $\phi\colon\Omega\to P$ be maps
such that $\phi\circ \psi$ is homotopic to the identity $\id_P$.
The locally trivial fibration over a bouquet
$p_\Omega=\phi^*(p_P)\colon E_\Omega\to\Omega$
has a global section if and only if it has a
section over every sphere of the bouquet.
If the fiber $M$ has non-abelian fundamental group
and is not homeomorphic to Klein bottle, then the space
of autohomeomorphisms $\Homeo(M)$ has simply connected identity
component~\cite{BH} and therefore any locally trivial fibration
over 2-sphere with fiber homeomorphic to $M$ has a section
(in fact, this fibration is trivial).
Hence, the fibration $p_\Omega$ has a section $s_\Omega$.
This section defines a lifting of the map $\phi\circ\psi\colon P\to P$
with respect to $p_P$.
Since $p_P$ is a Serre fibration and $\phi\circ\psi$
is homotopic to the identity, the identity mapping $\id_P$
has a lifting $s_P\colon P\to E_P$
with respect to $p_P$ which is simply a section of $p_P$.
\item[(c)]
Since $\pi_1(P)=\pi_2(P)=0$ and $\dim P=3$, then $P$ is homotopy equivalent
to a bouquet of 3-spheres $\Omega=\vee_{i=1}^m S^3_i$.
Let $\psi\colon P\to\Omega$ and $\phi\colon\Omega\to P$ be maps
such that $\phi\circ \psi$ is homotopic to the identity $\id_P$.
The locally trivial fibration over the bouquet
$p_\Omega=\phi^*(p_P)\colon E_\Omega\to\Omega$
has a global section if and only if it has a
section over every sphere of the bouquet.
Since the space of autohomeomorphisms of the Klein bottle
$\Homeo(K^2)$ has $\pi_2(\Homeo(K^2))=0$~\cite{BH},
any locally trivial fibration
over 3-sphere with fiber homeomorphic to $K^2$ has a section
(in fact, this fibration is trivial).
Hence, the fibration $p_\Omega$ has a section $s_\Omega$.
This section defines a lifting of the map $\phi\circ\psi\colon P\to P$
with respect to $p_P$.
Since $p_P$ is a Serre fibration and $\phi\circ\psi$
is homotopic to the identity, the identity mapping $\id_P$
has a lifting $s_P\colon P\to E_P$
with respect to $p_P$ which is simply a section of $p_P$.
\end{itemize}
\end{proof}

By the construction of $P$ the section $s_P$ defines an
$\varepsilon$-section for $p_L$. Therefore, $p_L$ has a section $s_L$.
Clearly, $s_L$ defines a lifting $T\colon L\to E$
of $\mu$ with respect to $p$.
Finally, we define compact singular filtration
$\mathcal F=\{(F_i,T_i)\}_{i=0}^2$ of $F$ as follows:
\[  F_0=F_1=\mu^{-1}\colon B\to L,\qquad F_2=F,
\qquad T_i=\id \text{ for } i=0  \]
and $T_1$ is defined fiberwise by
$T_1(x)=T|_{\mu^{-1}(x)}\colon \mu^{-1}(x)\to F(x)$.
The filtration $\mathcal F$ is approximately connected since for
$i=0,1$ any compactum $F_i(x)$ is $UV^1$.
And $\mathcal F$ is approximately $\infty$-connected since
every compactum $F(x)$ is an aspheric 2-manifold
(and therefore is approximately $n$-aspheric for all $n\ge 2$).

Now we can apply Theorem~\ref{thmapproxsections} to the mapping $F$
to find a single-valued continuous selection $s\colon B\to l_2$ of $F$.
Clearly, $s$ defines a section of the fibration $p$.
\end{proof}

The following Remark explains the appearence of the condition (c)
in Theorem~\ref{glob2sec}.

\begin{rem}
There exists a locally trivial fibration over 2-sphere with fibers
homeomorphic to Klein bottle having no global section.
\end{rem}

\section{Fibrations with 3-manifold fibers}\label{34}

Our strategy of construction the section here is the same
as in Section~\ref{33}.
The only difference is that instead of local hereditary asphericity
of fibers we will use local hereditary coconnected asphericity.

\begin{lem} \label{lemma3shrinking}
Let $F\colon X\to Y$ be
equi locally hereditarily coconnectedly aspheric,
lower coconnected, equi locally compact,
lower 2-continuous complete multivalued mapping
of $ANR$-space $X$ to Banach space $Y$.
Suppose that a compact submapping $\Psi\colon A\to Y$ of $F|_A$
is defined on compactum $A\subset X$ and admits continuous approximations.
Then for any $\varepsilon>0$ there exists a neighborhood $OA$ of $A$
and a compact submapping $\Psi'\colon OA\to Y$ of $F|_{OA}$
such that $\Gamma_{\Psi'}\subset O(\Gamma_\Psi,\varepsilon)$,
$\Psi'$ admits a compact approximately $\infty$-connected
filtration, and $\cal \Psi'<\varepsilon$.
\end{lem}

\begin{proof}
Fix a positive number $\varepsilon$.
Apply Lemma~\ref{L32} with $\alpha$ being equi local compactness
to get a positive number $\varepsilon_4<\varepsilon$.
Apply Lemma~\ref{L32} with $\alpha$ being equi local hereditary
coconnected asphericity to get a positive number $\varepsilon_3<\varepsilon_4$.
Apply Lemma~\ref{L32} with $\alpha$ being lower coconnectedness
to get a positive number $\varepsilon_2<\varepsilon_3/2$.
By Lemma~\ref{lempolcont} the mapping $F$ is lower polyhedrally 2-continuous.
Subsequently applying Lemma~\ref{L32} with $\alpha$ being
polyhedral $n$-continuity for $n=2,1,0$, we find positive numbers
$\varepsilon_1$, $\varepsilon_0$, and $\delta$ such that $\delta<\varepsilon_0<\varepsilon_1<\varepsilon_2$
and for every point $(x,y)\in O(\Gamma_{\Psi},\delta)$
the pair $(O(y,\varepsilon_1)\cap F(x),O(y,\varepsilon_2)\cap F(x))$
is polyhedrally 2-connected,
the pair $(O(y,\varepsilon_0)\cap F(x),O(y,\varepsilon_1)\cap F(x))$
is polyhedrally 1-connected, and
the intersection $O(y,\varepsilon_0)\cap F(x)$ is not empty.

Let $f\colon K\to Y$ be a continuous single-valued mapping
whose graph is contained in $O(\Gamma_{\Psi},\delta)$.
Let $f'\colon\mathcal OK\to Y$ be a continuous extension of the mapping $f$
over some neighborhood $\mathcal OK$ such that the graph of $f'$
is contained in $O(\Gamma_{\Psi},\delta)$.
Now we can define a polyhedrally connected filtration
$G_0\subset G_1\subset G_2\colon \mathcal OK\to Y$
of the mapping $F|_{\mathcal OK}$ by the equality
\[
G_i(x)=O(f'(x),\varepsilon_i)\cap F(x).
\]
Since the set $\cup_{x\in\mathcal OK}\{\{x\}\times O(f'(x),\varepsilon_i)\}$
is open in the product $\mathcal OK\times Y$ and the mapping $F$
is complete, then $G_i$ is also complete.
Clearly, the set $G_2^\Gamma(x)$ is contained in $O(\Gamma_\Psi,2\varepsilon_2)$.
Now, applying Lemma~\ref{L31} to the filtration $G_0\subset G_1\subset G_2$,
we obtain a compact approximately connected subfiltration
$F_0\subset F_1\subset F_2\colon \mathcal OK\to Y$.
By the choice of $\varepsilon_2$ we can find u.s.c. closed-valued mapping
$F_3$ containing $F_2$ such that the pair $F_2(x)\subset F_3(x)$
is coconnected in $F(x)$ for any $x\in X$.
By the choice of $\varepsilon_4$ we may assume that $F_3$ is compact.
Using Lemma~\ref{lemmacocon} and the choice of $\varepsilon_3$
we find a coconnectification $F_4$ of $F_2$ inside $F_3$.
Then $F_4$ is compact submapping of $F$ having approximately aspheric
point-images. Therefore, the filtration $F_0\subset F_1\subset F_4$
is approximately $\infty$-connected and we can put $\Psi'=F_4$.
\end{proof}

\begin{thm} \label{thm3approxsections}
Let $F\colon X\to Y$ be
equi locally hereditarily coconnectedly aspheric,
lower coconnected, equi locally compact,
lower 2-continuous complete multivalued mapping
of locally compact $ANR$-space $X$ to Banach space $Y$.
Suppose that a compact submapping $\Psi\colon A\to Y$ of $F|_A$
is defined on compactum $A\subset X$ and admits continuous approximations.
Then for any $\varepsilon>0$ there exists a neighborhood $OA$ of $A$
and a single-valued continuous selection $s\colon OA\to Y$
of $F|_{OA}$ such that $\Gamma_s\subset O(\Gamma_\Psi,\varepsilon)$.
\end{thm}

\begin{proof}
Consider a $G_\delta$-subset $G\subset X\times Y$ such that
all fibers of $F$ are closed in $G$ and fix open sets
$G_i\subset X\times Y$ such that $G=\cap^\infty_{i=1} G_i$.
Fix $\varepsilon>0$ such that $O(\Gamma_\Psi,\varepsilon)\subset G_1$.
By Lemma~\ref{lemma3shrinking} there is a neighborhood $U_1$
of $A$ in $X$ and a compact submapping $\Psi_1\colon U_1\to Y$ of $F|_{U_1}$
such that $\Gamma_{\Psi_1}\subset O(\Gamma_\Psi,\varepsilon)$,
$\Psi_1$ admits a compact approximately $\infty$-connected
filtration, and $\cal \Psi_1<\varepsilon$.
Since $X$ is locally compact and $A$ is compact, there exists a compact
neighborhood $OA$ of $A$ such that $OA\subset U_1$.
By Theorem~\ref{thmapprox} the mapping $\Psi_1|_{OA}$
admits continuous approximations.
Take $\varepsilon_1<\varepsilon$ such that the neighborhood
$\mathcal U_1=O(\Gamma_{\Psi_1}(OA),\varepsilon_1)$ lies in $O(\Gamma_\Psi,\varepsilon)$.
Clearly, $\mathcal U_1\subset G_1$.

Now by induction with the use of Lemma~\ref{lemma3shrinking},
we construct a sequence of neighborhoods
$U_1\supset U_2\supset U_3\supset\dots$ of the
compactum $OA$, a sequence of compact submappings
$\{\Psi_k\colon U_k\to Y\}^\infty_{k=1}$ of the mapping $F$,
and a sequence of neighborhoods
$\mathcal U_k=O(\Gamma_{\Psi_1}(OA),\varepsilon_k)$
such that for every $k\ge 2$ we have $\cal \Psi_k<\varepsilon_{k-1}/2<\varepsilon/2^k$,
and $\mathcal U_k(OA)$ is contained in $\mathcal U_{k-1}(OA)\cap G_k$.
It is not difficult to choose the neighborhood $\mathcal U_k$
of the graph $\Gamma_{\Psi_k}$ in such a way that for every point $x\in U_k$
the set $\mathcal U_k(x)$ has diameter less than $\frac{3}{2^k}$.

Then for any $m\ge k\ge 1$ and for any point $x\in OA$ we have
$\Psi_m(x)\subset O(\Psi_k(x),\frac{3}{2^k})$; this implies the fact
that the sequence $\{\Psi_k|_{OA}\}_{k=1}^\infty$ is a Cauchy sequence.
Since $Y$ is complete, there exists a limit $s\colon OA\to Y$ of this sequence.
The mapping $s$ is single-valued by
the condition $\cal\Psi_k<\frac{1}{2^k}$ and is upper
semicontinuous (and, therefore, is continuous) by the upper
semicontinuity of all the mappings $\Psi_k$.
Clearly, for any $x\in OA$ the point $s(x)$ lies in $G(x)$
and is a limit point of the set $F(x)$.
Since $F(x)$ is closed in $G(x)$, then $s(x)\in F(x)$,
i.e. $s$ is a selection of the mapping $F$.
\end{proof}

\begin{thm}
Let $p\colon E\to B$ be a topologically regular mapping
of compacta with fibers homeomorphic to a 3-dimensional manifold.
If $B\in ANR$, then any section of $p$ over closed subset $A\subset B$
can be extended to a section of $p$ over some neighborhood of $A$.
\end{thm}

\begin{proof}
Let $s\colon A\to E$ be a section of $p$ over $A$.
Embed $E$ into Hilbert space $l_2$ and consider a multivalued mapping
$F\colon B\to l_2$ defined as follows:
$$
   F(b)=\cases s(b), &\text{if \;$b\in A$}\\
               p^{-1}(b),   &\text{if \;$x\in B\setminus A$.}\endcases
$$
Since every fiber $p^{-1}(b)$ is compact, the mapping $F$ is complete.
It follows from Lemma~\ref{lemmafibrprop} that the mapping $F$ is
equi locally hereditarily coconnectedly aspheric,
lower coconnected, equi locally compact, and lower $2$-continuous.
We can apply Theorem~\ref{thm3approxsections} to the mapping $F$
and its submapping $s$ to find a single-valued continuous selection
$\widetilde s\colon OA\to l_2$ of $F|_{OA}$.
By definition of $F$, we have $\widetilde s|_{A}=F|_{A}=s$.
Clearly, $\widetilde s$ defines a section of the fibration $p$
over $OA$ extending $s$.
\end{proof}

\begin{thm}
Let $p\colon E\to B$ be topologically regular mapping of compacta
with fibers homeomorphic to some compact connected 3-dimensional manifold $M$.
If $B$ is ANR-space, then $p$ admits a global section
if either of the following conditions hold:
\begin{itemize}
\item[(a)] $\pi_1(M)$ is abelian, $M$ is aspheric, and $H^2(B;\pi_1(F_b))=0$
\item[(b)] $M$ is closed hyperbolic 3-manifold and $\pi_1(B)=0$.
\item[(c)] $M$ is closed, irreducible, sufficiently large, contains
no embedded $\mathbb RP^2$ having a trivial normal bundle, and
$\pi_1(B)=\pi_2(B)=0$.
\end{itemize}
\end{thm}

\begin{proof}
Embed $E$ into Hilbert space $l_2$ and consider a multivalued mapping
$F\colon B\to l_2$ defined as $F=p^{-1}$.
Since every fiber $p^{-1}(b)$ is compact, the mapping $F$ is complete.
It follows from Lemma~\ref{lemmafibrprop} that the mapping $F$ is
equi locally hereditarily coconnectedly aspheric,
lower coconnected, equi locally compact, and lower $2$-continuous.

Now we show that $F$ admits a compact singular
approximately $\infty$-connected filtration.
In cases (a) and (b) there exists $UV^1$-mapping $\mu$
of Menger 2-dimensional manifold $L$
onto $B$~\cite{Dr1}. Note that $\pi_1(L)=0$ if $\pi_1(B)=0$.
In case (c) we consider $UV^2$-mapping $\mu$
of Menger 3-dimensional manifold $L$ onto $B$~\cite{Dr1};
note that $\pi_1(L)=\pi_2(L)=0$ if $\pi_1(B)=\pi_2(B)=0$.
Since $\dim L<\infty$, the induced fibration $p_L=\mu^*(p)\colon E_L\to L$
is locally trivial~\cite{H}.
By Proposition~\ref{profibr} there is $\varepsilon>0$ such that an
existence of $\varepsilon$-section for $p_L$ implies an
existence of a section for $p_L$.
In cases (a) and (b) by Proposition~\ref{propoly} there exist a 2-dimensional
finite polyhedron $P$ and continuous maps $g\colon L\to P$
and $h\colon P\to L$ such that $h\circ g$ is $\varepsilon$-close
to the identity (we assume $\pi_1(P)=0$ in case $\pi_1(B)=0$).
In case (c) by Proposition~\ref{pro3poly} there exist a 3-dimensional
finite polyhedron $P$ with $\pi_1(P)=\pi_2(P)=0$
and continuous maps $g\colon L\to P$ and $h\colon P\to L$
such that $h\circ g$ is $\varepsilon$-close to the identity.

Consider a locally trivial fibration $p_P=h^*(p_L)\colon E_P\to P$.

{\bf Claim.} \ {\it The fibration $p_P$ has a section $s_P$.}

\begin{proof}
\begin{itemize}
\item[(a)]
If $\pi_1(M)$ is abelian and $H^2(B;\pi_1(F_b))=0=H^2(P;\pi_1(F_b))$,
the fibration $p_P$ has a section $s_P$~\cite{W}.

\item[(b)]
Since $\pi_1(P)=0$ and $\dim P=2$, then $P$ is homotopy equivalent
to a bouquet of 2-spheres $\Omega=\vee_{i=1}^m S^2_i$.
Let $\psi\colon P\to\Omega$ and $\phi\colon\Omega\to P$ be maps
such that $\phi\circ \psi$ is homotopic to the identity $\id_P$.
The locally trivial fibration over a bouquet
$p_\Omega=\phi^*(p_P)\colon E_\Omega\to\Omega$
has a global section if and only if it has a
section over every sphere of the bouquet.
For a closed hyperbolic 3-manifold $M$ the space
of autohomeomorphisms $\Homeo(M)$ has simply connected identity
component~\cite{G} and therefore any locally trivial fibration
over 2-sphere with fiber homeomorphic to $M$ has a section
(in fact, this fibration is trivial).
Hence, the fibration $p_\Omega$ has a section $s_\Omega$.
This section defines a lifting of the map $\phi\circ\psi\colon P\to P$
with respect to $p_P$.
Since $p_P$ is a Serre fibration and $\phi\circ\psi$
is homotopic to the identity, the identity mapping $\id_P$
has a lifting $s_P\colon P\to E_P$
with respect to $p_P$ which is simply a section of $p_P$.

\item[(c)]
Since $\pi_1(P)=\pi_2(P)=0$ and $\dim P=3$, then $P$ is homotopy equivalent
to a bouquet of 3-spheres $\Omega=\vee_{i=1}^m S^3_i$.
Let $\psi\colon P\to\Omega$ and $\phi\colon\Omega\to P$ be maps
such that $\phi\circ \psi$ is homotopic to the identity $\id_P$.
The locally trivial fibration over the bouquet
$p_\Omega=\phi^*(p_P)\colon E_\Omega\to\Omega$
has a global section if and only if it has a
section over every sphere of the bouquet.
Since the space of autohomeomorphisms $\Homeo(M)$
in this case has $\pi_2(\Homeo(M))=0$~\cite{Hat},
any locally trivial fibration
over 3-sphere with fiber homeomorphic to $M$ has a section
(in fact, this fibration is trivial).
Hence, the fibration $p_\Omega$ has a section $s_\Omega$.
This section defines a lifting of the map $\phi\circ\psi\colon P\to P$
with respect to $p_P$.
Since $p_P$ is a Serre fibration and $\phi\circ\psi$
is homotopic to the identity, the identity mapping $\id_P$
has a lifting $s_P\colon P\to E_P$
with respect to $p_P$ which is simply a section of $p_P$.
\end{itemize}
\end{proof}

By the construction of $P$ the section $s_P$ defines an
$\varepsilon$-section for $p_L$. Therefore, $p_L$ has a section $s_L$.
Clearly, $s_L$ defines a lifting $T\colon L\to E$
of $\mu$ with respect to $p$.
Finally, we define compact singular filtration
$\mathcal F=\{(F_i,T_i)\}_{i=0}^2$ of $F$ as follows:
\[  F_0=F_1=\mu^{-1}\colon B\to L,\qquad F_2=F,
\qquad T_i=\id \text{ for } i=0  \]
and $T_1$ is defined fiberwise by
$T_1(x)=T|_{\mu^{-1}(x)}\colon \mu^{-1}(x)\to F(x)$.
The filtration $\mathcal F$ is approximately connected since for
$i=0,1$ any compactum $F_i(x)$ is $UV^1$.
And $\mathcal F$ is approximately $\infty$-connected since
every compactum $F(x)$ is an aspheric 2-manifold
(and therefore is approximately $n$-aspheric for all $n\ge 2$).

Now we can apply Theorem~\ref{thm3approxsections} to the mapping $F$
to find a single-valued continuous selection $s\colon B\to l_2$ of $F$.
Clearly, $s$ defines a section of the fibration $p$.
\end{proof}

\section{Acknowledgements}

Authors wish to express their sincere thanks to
P.~Akhmetiev, R.J.~Daverman, B.~Hajduk and T.~Yagasaki
for helpful discussions during the development of this work.

\end{document}